\relax 
\global \edef \factsub {{0.1}}
\global \edef \factcolony {{0.3}}
\global \edef \factfewvars {{0.4}}
\global \edef \factassocses {{0.5}}
\global \edef \sectnewfam {{1}}
\global \edef \qsectnewfam {{2}}
\global \edef \thmdbl {{1.1}}
\global \edef \Qfivesix {{7}}
\global \edef \propIrr {{2.2}}
\global \edef \qKEFcolon {{9}}
\global \edef \qV {{9}}
\global \edef \propIrrr {{2.3}}
\global \edef \qcalco {{11}}
\global \edef \qKEFVcolon {{14}}
\global \edef \propIrrrr {{2.4}}
\global \edef \propIrrrrr {{2.5}}
\global \edef \thmthree {{2.6}}
\global \edef \thmtwo {{2.7}}
\magnification=\magstep1
\baselineskip=15pt
\overfullrule=0pt

\font\footfont=cmr8
\font\footitfont=cmti8
\font\sc=cmcsc10 at 12pt
\font\large=cmbx12 at 14pt

\font\tenmsb=msbm10
\font\sevenmsb=msbm7
\font\fivemsb=msbm5
\newfam\msbfam
\textfont\msbfam=\tenmsb
\scriptfont\msbfam=\sevenmsb
\scriptscriptfont\msbfam=\fivemsb
\def\hexnumber#1{\ifcase#1 0\or1\or2\or3\or4\or5\or6\or7\or8\or9\or
	A\or B\or C\or D\or E\or F\fi}

\mathchardef\subsetneq="2\hexnumber\msbfam28


\def\Ass{\hbox{Ass}\hskip -0.0em}

\def\qed{\hbox{\quad \vrule width 1.6mm height 1.6mm depth 0mm}\vskip 1ex}
\def\eqed{\hbox{\quad \vrule width 1.6mm height 1.6mm depth 0mm}}

\def\D{D}
\def\C{C}

\def\newembprno{$160 n - 301 + 16 d + n(n-1)
+ 31(d^{2^1}+ \cdots + d^{2^{n-3}})
+(n-1)d^{2^1} +(n-2) d^{2^2} + \cdots + 3d^{2^{n-3}}
+18 d^{2^{n-2}}$}

\newcount\tableco \tableco=0
\def\tabledo{\global\advance\tableco by 1{\the\tableco}}
\def\tabel#1{{\global\edef#1{\the\tableco}}}

\def\label#1{{\global\edef#1{\the\sectno.\the\thmno}}\ignorespaces}
\def\pagelabel#1{\immediate\write\isauxout{\noexpand\global\noexpand\edef\noexpand#1{{\the\pageno}}}{\global\edef#1{\the\pageno}}\ignorespaces}

\def\isnameuse#1{\csname #1\endcsname}

\def\issecond#1#2{#2}
\def\isifundefined#1#2#3{
	\expandafter\ifx\csname #1\endcsname\relax #2
	\else #3 \fi}
\def\pageref#1{\isifundefined{is#1}
	{{\bf ??}\message{Reference `#1' on page [\number\count0] undefined}}
	{\edef\istempa{\isnameuse{is#1}}\expandafter\issecond\istempa\relax}}

\def\today{\ifcase\month\or January\or February\or March\or
  April\or May\or June\or July\or August\or September\or
  October\or November\or December\fi
  \space\number\day, \number\year}

\newcount\sectno \sectno=0
\newcount\thmno \thmno=0
\def \section#1{\vskip 1.2truecm
	\global\advance\sectno by 1 \global\thmno=0
	\noindent{\bf \the\sectno. #1} \vskip 0.6truecm}
\def \thmline#1{\vskip 15pt
	\global\advance\thmno by 1
	\noindent{\bf #1\ \the\sectno.\the\thmno:}\ \ %
	\bgroup \advance\baselineskip by -1pt \it
	\abovedisplayskip =4pt
	\belowdisplayskip =3pt
	\parskip=0pt
	}

\def \thm{\thmline{Theorem}}

\def \endb{\egroup \vskip 1.4ex}

\def \fact{\advance\thmno by 1\item{\bf\the\sectno.\the\thmno:}}
\def \prop{\thmline{Proposition}}

\def \proofof#1{\medskip\noindent {\sl Proof of #1:\ }}


\def\label#1{\unskip\immediate\write\isauxout{\noexpand\global\noexpand\edef\noexpand#1{{\the\sectno.\the\thmno}}}
    {\global\edef#1{\the\sectno.\the\thmno}}\unskip\ignorespaces}
\def\sectlabel#1{\immediate\write\isauxout{\noexpand\global\noexpand\edef\noexpand#1{{\the\sectno}}}
    {\global\edef#1{\the\sectno}}}

\def\tabel#1{\hbox to 0pt{\hskip -6em\string#1}\unskip\immediate\write\isauxout{\noexpand\global\noexpand\edef\noexpand#1{{\the\tableco}}}\unskip}

\newwrite\isauxout
\openin1\jobname.aux
\ifeof1\message{No file \jobname.aux}
       \else\closein1\relax\input\jobname.aux
       \fi
\immediate\openout\isauxout=\jobname.aux
\immediate\write\isauxout{\relax}

\ %
\vskip 3ex
\centerline{\large A new family of ideals}
\centerline{\large with the doubly exponential ideal membership property}

\vskip 4ex
\centerline{\sc Irena Swanson}
\centerline{\sc \today}
\unskip\footnote{ }{{\footitfont 1991 Mathematics Subject Classification.}
13C13, 13P05}
\unskip\footnote{ }{{\footitfont Key words and phrases.}
\footfont Primary decomposition, Mayr-Meyer, membership problem, complexity of ideals.}

\vskip 0.5cm

Mayr and Meyer [MM] found ideals with the doubly exponential ideal membership property.
Further investigations of the doubly exponential properties
of these ideals can be found in Bayer and Stillman [BS], Demazure [D], and Koh [K].
Following a question of Bayer, Huneke and Stillman,
the author has investigated in [S1, S2, S3] the properties
of the primary decompositions and the associated primes of these Mayr-Meyer ideals.
In that investigation a new family of ideals arose.
This new family $K(n,d)$ is presented and analyzed in this paper.
It is proved in this paper
that this new family also satisfies the doubly exponential ideal membership
property.
Thus the question of Bayer, Huneke and Stillman also pertains to this family.
The main part of the paper is finding the associated primes of this family.

This new family is ideally suited -- perhaps by construction --
to recursively constructing a set containing all of its associated primes.
The recursively built set most likely contains prime ideals
which are not necessarily associated.
The elimination of redundancies is not attempted here.

The main tool used below for finding the associated primes of the Mayr-Meyer
ideals are various short exact sequences:
the associated primes of the middle module
is contained in the union of the associated primes of the two other modules.
Theorems~\thmthree\ and~\thmtwo\ give a set of prime ideals obtained in this way.
The total number of possibly embedded primes of $K(n,d)$ found in this paper
is \newembprno.
This number is doubly exponential in $n$.
It is not known whether the set of associated primes of $K(n,d)$
is indeed doubly exponential.
These same prime ideals,
after adding $6$ more variables,
are possibly also associated to the Mayr-Meyer ideals
(see [S3]).

\vskip 4ex

{\bf Acknowledgement.}
I thank Craig Huneke for suggesting this problem
all for all the conversations and enthusiasm for this research.
I also thank the NSF for partial support on grants DMS-0073140 and DMS-9970566.
Furthermore,
I am indebted to the symbolic computer algebra packages Macaulay2 and Singular
for verifications of the proofs below
for low values of $n$ and $d$.

The following summarizes the elementary facts about primary decompositions
used in this paper (as well as in [S1, S2, S3]):

\vskip2ex\noindent{\bf Facts:}
\bgroup\parindent=3em
\fact
\label{\factsub}
For any ideals $I, I'$ and $I''$ with $I \subseteq I''$,
$(I + I') \cap I'' = I + I' \cap I''$.

\fact
For any ideal $I$ and element $x$,
$(x) \cap I = x (I : x)$.

\fact
\label{\factcolony}
For any ideals $I$ and $I'$,
and any element $x$,
$(I + xI') : x = (I : x) + I'$.

\fact
\label{\factfewvars}
Let $x_1, \ldots, x_n$ be variables over a ring $R$.
Let $S = R[x_1, \ldots, x_n]$.
For any $f_1 \in R$,
$f_2 \in R[x_1]$,
$\ldots$, $f_n \in R[x_1, \ldots, x_{n-1}]$,
let $L$ be the ideal $(x_1 - f_1, \ldots, x_n - f_n)S$ in $S$.
Then an ideal $I$ in $R$ is primary (respectively, prime)
if and only if $IS + L$ is a primary (respectively, prime) in $S$.
Furthermore,
$\cap_i q_i = I$ is a primary decomposition of $I$
if and only if $\cap_i (q_iS + L)$ is a primary decomposition of $IS + L$.

\fact
\label{\factassocses}
Let $I$ be an ideal in a ring $R$.
Then for any $x \in R$,
$\Ass\left({R \over I}\right) \subseteq
\Ass\left({R \over I : x}\right) \cup
\Ass\left({R \over I+ (x)}\right)$,
and every associated prime of
${R \over I : x}$ is an associated prime of ${R \over I}$.
(Use the short exact sequence
$0 \longrightarrow
{R \over I : x} \longrightarrow
{R \over I} \longrightarrow
{R \over I + (x)} \longrightarrow 0$.)

\egroup

\section{The family, and its doubly exponential membership property}
\sectlabel{\sectnewfam}
\pagelabel{\qsectnewfam}

We will define a new two-parameter family of ideals,
and prove that this family satisfies the doubly exponential ideal membership property.
The doubly exponential ideal membership property was first addressed in Hermann's paper~[H]:
if $I$ is an ideal in an $n$-dimensional polynomial ring
over the field of rational numbers,
and $I$ is generated by polynomials $f_1, \ldots, f_k$ of degree at most $d$,
then it is possible to write $f = \sum r_i f_i$
such that each $r_i$ has degree at most $\deg f + (kd)^{(2^n)}$.
Thus Hermann proved that every ideal there satisfies a doubly exponential
ideal membership property.
Fortunately for computational purposes,
most (known) ideals satisfy an ideal membership property of much smaller complexity.
It was not until 1982 that an ideal with a doubly exponential ideal membership
property was indeed found by Mayr and Meyer (see [MM]).
Further analyses and modifications were obtained by
Bayer and Stillman [BS], Demazure [D], and Koh [K].

Here are the definitions of the new family:
let $n, d \ge 2$ be integers,
and $k$ a field.
Let 
$s, f, s_{r+1}, f_{r+1}$,
$b_{r1}, b_{r2}, b_{r3}, b_{r4}$,
$c_{r1}, c_{r2}, c_{r3}, c_{r4}$
be variables over $k$.
Set
$$
R = k[s_2, \ldots, s_n, f_2, \ldots, f_n,
b_{0i}, \ldots, b_{n-1,i}, c_{1i}, \ldots, c_{n-1,i}, | i = 1, \ldots, 4].
$$
Thus $R$ is a polynomial ring of dimension $10n-6$.
The ideal $K_l(n,d)$ in $R$ (``l" for long, a shortened version appearing later)
is generated by the following generators:
$$
\eqalignno{
& G_{01} = b_{01} b_{03}^d - b_{04} b_{02}^d, \cr
& G_{1i} = c_{1i} \left(b_{02}-b_{1i} b_{03}\right), i = 1, \ldots, 4, \cr
& G_{1,4+i} = c_{1i} \left(b_{01}-b_{1i}^d b_{04}\right), i = 1, \ldots, 4, \cr
& G_{1ij} = c_{1i} c_{1j} \left(b_{1i} - b_{1j}\right), 1 \le i < j \le 4, \cr
& G_{21} = b_{04}^d c_{11} - b_{01}^d c_{12}, \cr
& G_{22} = b_{04}^d c_{14} - b_{01}^d c_{13},\cr
& G_{23} = b_{01}^d (c_{12} -c_{13}), \cr
& G_{24} = b_{04}^d (c_{12} b_{11} -c_{13} b_{14}), \cr
& G_{2,4+i} = b_{04}^d c_{12} c_{2i} \left( b_{12}-b_{2i} b_{13} \right),
i = 1, \ldots, 4. \cr
&G_{20} = s_2 - c_{11}, \cr
&G_{2,-1} = f_2 - c_{14}, \cr
&G_{r0} = s_r - s_{r-1} c_{r-1,1}, r = 3, \ldots, n-1, \cr
&G_{r,-1} = f_r - s_{r-1} c_{r-1,4}, r = 3, \ldots, n-1, \cr
&G_{r1} = b_{01}^d \left(f_{r-1} c_{r-1,1} - s_{r-1} c_{r-1,2} \right), r = 3, \ldots, n, \cr
&G_{r2} = b_{01}^d \left(f_{r-1} c_{r-1,4} - s_{r-1} c_{r-1,3} \right), r = 3, \ldots, n, \cr
&G_{r3} = b_{01}^d s_{r-1} \left(c_{r-1,3} - c_{r-1,2} \right), r = 3, \ldots, n, \cr
&G_{r4} = b_{01}^d f_{r-1} \left(c_{r-1,2} b_{r-1,1} - c_{r-1,3} b_{r-1,4} \right), r = 3, \ldots, n, \cr
&G_{r,4+i} = b_{01}^d f_{r-1} c_{r-1,2} c_{ri} \left(b_{r-1,2} - b_{ri} b_{r-1,3} \right),
r = 3, \ldots, n-1, \cr
&G_{n0} = s_n - s_{n-1} c_{n-1,1} b_{01}^d, \cr
&G_{n,-1} = f_n - s_{n-1} c_{n-1,4} b_{01}^d, \cr
&G_{n5} = b_{01}^d f_{n-1} c_{n-1,2}  \left(b_{n-1,2} - b_{n-1,3} \right). \cr
}
$$

Note that the maximal degree of a generator of $K_l(n,d)$ is $d + 5$.

\thm
\label{\thmdbl}
With notation as above,
the element $s_n - f_n$ of $R$ lies in $K_l(n,d)$,
has degree $1$,
but when written as an $R$-linear combination of the given generators of $K_l(n,d)$,
the degree of at least one coefficient is doubly exponential in $n$.

In other words,
the family $K_l(n,d)$ satisfies the doubly exponential ideal membership property.
\endb

We postpone the proof in favor of first introducing more notation.

Under the evaluation mapping
$$
\eqalignno{
&s_r \mapsto c_{11} \cdots c_{r-1,1},
f_r \mapsto c_{11} \cdots c_{r-2,1} c_{r-1,4},
r = 2, \ldots, n-1, \cr
& s_n \mapsto c_{11} \cdots c_{n-1,1} b_{01}^d,
f_n \mapsto c_{11} \cdots c_{n-2,1} c_{n-1,4} b_{01}^d, \cr
}
$$
the image of $K_l(n,d)$ is an ideal $K(n,d)$ in $R$
(or actually in the polynomial subring of $R$ obtained by omitting the $s_r$, $f_r$)
generated by the following elements:
first one level 0 generator:
$$
\eqalignno{
g_{01} &= b_{01} b_{03}^d - b_{04} b_{02}^d, \cr
}
$$
then fourteen level 1 generators:
$$
\eqalignno{
& g_{1i} = c_{1i} \left(b_{02}-b_{1i} b_{03}\right), i = 1, \ldots, 4, \cr
& g_{1,4+i} = c_{1i} \left(b_{01}-b_{1i}^d b_{04}\right), i = 1, \ldots, 4, \cr
& g_{1ij} = c_{1i} c_{1j} \left(b_{1i} - b_{1j}\right), 1 \le i < j \le 4, \cr
}
$$
then eight level 2 generators:
$$
\eqalignno{
& g_{21} = b_{04}^d c_{11} - b_{01}^d c_{12}, \cr
& g_{22} = b_{04}^d c_{14} - b_{01}^d c_{13},\cr
& g_{23} = b_{01}^d (c_{12} -c_{13}), \cr
& g_{24} = b_{04}^d (c_{12} b_{11} -c_{13} b_{14}), \cr
& g_{2,4+i} = b_{04}^d c_{12} c_{2i}
\left( b_{12}-b_{2i} b_{13} \right),
i = 1, \ldots, 4. \cr
}
$$
When $n = 2$,
the last four generators are replaced by only one:
$$
g_{25} = b_{04}^d c_{12} c_{2i}\left( b_{12}-b_{13} \right).
$$
Next,
the first four level $r$ generators, $r = 3, \ldots, n$, are:
$$
\eqalignno{
&g_{r1} = b_{01}^d c_{11} \cdots c_{r-3,1} \left(
c_{r-2,4} c_{r-1,1} - c_{r-2,1} c_{r-1,2}
\right), \cr
& g_{r2} = b_{01}^d c_{11} \cdots c_{r-3,1} \left(
c_{r-2,4} c_{r-1,4} - c_{r-2,1} c_{r-1,3}
\right), \cr
& g_{r3} = b_{01}^d c_{11} \cdots c_{r-2,1}
\left( c_{r-1,3}-c_{r-1,2} \right), \cr
& g_{r4} = b_{01}^d c_{11} \cdots c_{r-3,1} c_{r-2,4}
\left( c_{r-1,2} b_{r-1,1}-c_{r-1,3} b_{r-1,4} \right), \cr
}
$$
the last four level $r$ generators, $r= 2, ..., n-1$, are:
$$
g_{r,4+i} = b_{01}^d c_{11} \cdots c_{r-3,1} c_{r-2,4} c_{r-1,2} c_{ri}
\left( b_{r-1,2}-b_{ri} b_{r-1,3} \right),
i = 1, \ldots, 4,
$$
and the last level $n$ generator is:
$$
g_{n5} = b_{01}^d c_{11} \cdots c_{n-3,1} c_{n-2,4} c_{n-1,2}
\left( b_{n-1,2}-b_{n-1,3} \right).
$$
As was computed in Section~5 of [S3],
$$
(J(n,d) : sc_{02}) + (c_{02},f) = K(n,d) + (s,f,c_{01}, c_{02}, c_{03}, c_{04}),
$$
where $J(n,d)$ is the original Mayr-Meyer ideal
defined in the polynomial ring $R[s,f,$
$c_{01}, c_{02}, c_{03}, c_{04}]$.

With this we are ready to prove Theorem \thmdbl:

\proofof{\thmdbl}
By the definitions,
the element $s_n - f_n$ is in $K_l(n,d)$ if and only if
the element $b_{01}^d c_{11} \cdots c_{n-2,1} (c_{n-1,1} - c_{n-1,4})$ is in $K(n,d)$.
As
$b_{01}^d c_{11} \cdots c_{n-2,1} (c_{n-1,1} - c_{n-1,4})$ and the generators of $K(n,d)$
do not involve any variables $s,f,c_{01}, c_{02}, c_{03}, c_{04}$,
then
$b_{01}^d c_{11} \cdots c_{n-2,1} (c_{n-1,1} - c_{n-1,4})$ is in $K(n,d)$
if and only if
$b_{01}^d c_{11} \cdots c_{n-2,1} (c_{n-1,1} - c_{n-1,4})$ is in
$$
K(n,d) + (s,f,c_{01}, c_{02}, c_{03}, c_{04})
= (J(n,d) : sc_{02}) + (c_{02},f),
$$
and then analogously this holds if and only if
$b_{01}^d c_{11} \cdots c_{n-2,1} (c_{n-1,1} - c_{n-1,4})$ is in $J(n,d) : sc_{02}$.
But this is indeed the case by Corollary~5.4 in [S3].

Thus we can write $s_n - f_n = \sum_{ri} A_{ri} G_{ri} + \sum_{ij} A_{1ij} G_{1ij}$,
for some $A_{ri}, A_{1ij} \in R$.
Let $N(n,d)$ be the maximum degree of the $A_{ri}, A_{1ij}$
with $i \ge 1$.

Under the evaluation map as in the discussion above this says that
$$
b_{01}^d c_{11} \cdots c_{n-2,1} (c_{n-1,1} - c_{n-1,4})
= \sum_{r;i\ge 1} a_{ri} g_{ri} + \sum_{ij} a_{1ij} g_{1ij},
$$
for some $a_{ri}, a_{1ij} \in R$ which do not involve any of the $s_r$, $f_r$.
Note that the degrees of the $a_{ri}$ and the $a_{1ij}$
are at most $N(n,d)\cdot(n-1+d)$.
Multiplying through by $sc_{02}$ gives
$$
b_{01}^d sc_{02} c_{11} \cdots c_{n-2,1} (c_{n-1,1} - c_{n-1,4})
= \sum_{r;i\ge 1} a_{ri} sc_{02} g_{ri} + \sum_{ij} a_{1ij} sc_{02} g_{1ij}.
$$
By Corollary~5.4 in [S3],
each of $s(c_{01} - b_{02}^d c_{02})$,
$sc_{02} g_{ri}$ and $sc_{02} g_{1ij}$ can be rewritten as
a linear combination of the generators of the Mayr-Meyer ideal
with coefficients of those generators having degrees at most $2d+1$.
In particular,
$$
\eqalignno{
sc_{01} c_{11} &\cdots c_{n-2,1} (c_{n-1,1} - c_{n-1,4})
=
s(c_{01} - b_{01}^d c_{02}) c_{11} \cdots c_{n-2,1} (c_{n-1,1} - c_{n-1,4}) \cr
&\hskip2em
+ b_{01}^d sc_{02} c_{11} \cdots c_{n-2,1} (c_{n-1,1} - c_{n-1,4}) \cr
&=
s(c_{01} - b_{01}^d c_{02}) c_{11} \cdots c_{n-2,1} (c_{n-1,1} - c_{n-1,4})
+\sum_{r;i\ge 1} a_{ri} sc_{02} g_{ri} + \sum_{ij} a_{1ij} sc_{02} g_{1ij} \cr
}
$$
can be written as a linear combination of the generators of the Mayr-Meyer ideal
whose coefficients have degrees at most 
$2d + 1 + \max\{n, \deg a_{ri}, \deg a_{1ij}\} \le 2d + 1 + \max\{n, N(n,d)(n-1+d)\}$.
By the work of Mayr and Meyer,
this maximum is in fact doubly exponential in $n$,
so that also $N(n,d)$ is doubly exponential in $n$.
\qed

\section{The associated primes of this family}

Bayer, Huneke and Stillman asked whether
the doubly exponential behavior of the Mayr-Meyer ideals
is reflected in some way in their associated primes.
The same question can be applied also to the family of ideals $K_l(n,d)$
defined in this paper.

As far as the primary decomposition of $K_l(n,d)$ is concerned,
by Fact~\factfewvars\ it suffices to find a primary decomposition of $K(n,d)$,
and then add some variables to the primary components.
The same holds for the associated primes of $K_l(n,d)$.
As $K(n,d)$ is notationally simpler and more relevant than $K_l(n,d)$ (for these purposes),
in the sequel we only work with $K(n,d)$.

For simplicity of notation we will assume in this section
that $k$ is an algebraically closed field
whose characteristic is relatively prime to $d$.

Sometimes we will write the ideal $K(n,d)$
also $K(n,d; c_{r'i}, b_{ri} | r' = 1, \ldots, n-1; r = 0, \ldots, n; i = 1, \ldots, 4)$,
to emphasize the defining variables.

\def\KEF{K(n,d)}
\def\KK{L}
\def\FF{M}
\def\EE{N}
\def\KKpone{\KK_1}
\def\EEpone{\EE_1}
\def\FFpone{\FF_1}
\def\KEpone{\KKpone + \EEpone}

Let $\FF$ be the ideal generated by all the level 0 and level 1
generators of $K(n,d)$,
let $\EE$ denote the subideal generated by all the level two generators $g_{2j}$,
and let $\KK$ denote the subideal of $K(n,d)$ generated
by all the $g_{rj}$ with $r \ge 3$.
Thus $K(n,d) = M + N + L$.

These new ideals $\FF, \EE$ and $\KK$ will sometimes be specified with
the variables and degrees
$(n,d; c_{r'i}, b_{ri} | r' = 1, \ldots, n-1; r = 0, \ldots, n; i = 1, \ldots, 4)$
attached to them.

Then define $\FFpone, \EEpone$ and $\KKpone$
as the corresponding level ideals of
$$
K(n-1,d^2) = K(n-1,d^2; c_{r'i}, b_{ri} |
r = 1, \ldots, n-1;
r' = 2, \ldots, n-1;
i = 1, \ldots, 4).
$$
Then $K(n-1,d^2)$
equals $\FFpone + \EEpone + \KKpone$.

We will prove that finding the set of associated primes of $K(n,d)$
reduces to finding the set of associated primes of
$(c_{11}, c_{12}, c_{13}, c_{14},$
$b_{01}, b_{02}, b_{03}, b_{04}) + K(n-1,d^2)$,
and by Fact~\factfewvars,
this reduces to finding the set of associated primes of $K(n-1,d^2)$.

Thus the associated prime ideals of $K(n,d)$ follow a recursion pattern.

For notational purposes we also define the following ideals in $R$:
$$
\C_r = (c_{r1}, c_{r2}, c_{r3}, c_{r4}), r = 1, \ldots, n-1.
$$

We now start the analysis of the associated primes of $K(n,d)$.
The main tool is Fact~\factassocses.

By Fact~\factassocses,
$\Ass \left( {R \over \KEF}\right) \subseteq
\Ass \left( {R \over \KEF+ (b_{04}^d)}\right)
\cup \Ass \left( {R \over \KEF: b_{04}^d}\right)$.
But $\KK \subseteq \C_1b_{01}^d\subseteq \C_1b_{04}^d +\FF$,
so that
$$
\eqalignno{
\KEF + (b_{04}^d) &= \EE+ \FF +(b_{04}^d) \cr
&=(b_{04}^d,
b_{01}b_{03}^d-b_{04}b_{02}^d) + c_{1i}(
b_{02}-b_{1i}b_{03}, c_{1j}(b_{1i}-b_{1j}),
b_{01}-b_{1i}^db_{04}) \cr
&= \bigcap_{\Lambda} \bigl(
(b_{04}^d,
b_{01}b_{03}^d-b_{04}b_{02}^d)
+ (c_{1i} | i \not \in \Lambda) \cr
&\hskip2em+ (b_{02}-b_{1i}b_{03},
b_{1i}-b_{1j},
b_{01}-b_{1i}^db_{04} | i, j \in \Lambda
) \bigr) \cr
&= \bigcap_{\Lambda\not =\emptyset} \bigl(
(b_{04}^d)
+ (c_{1i} | i \not \in \Lambda) + (b_{02}-b_{1i}b_{03},
b_{1i}-b_{1j},
b_{01}-b_{1i}^db_{04} | i, j \in \Lambda
) \bigr) \cr
&\hskip1em \cap \bigl( \C_1 +
(b_{01}^d,
b_{04}^d,
b_{01}b_{03}^d-b_{04}b_{02}^d) \bigr) \cap \bigl( \C_1 +
(b_{04}^d,
b_{03}^{d^2},
b_{01}b_{03}^d-b_{04}b_{02}^d)
\bigr) .\cr
}
$$
Thus the associated primes of $\KEF + (b_{04}^d)$ are
$$
\eqalignno{
Q_{1\Lambda} &=
\left(b_{01},b_{04}\right)
+ (c_{1i} | i \not \in \Lambda) + \left(b_{02}-b_{1i}b_{03},
b_{1i}-b_{1j} | i, j \in \Lambda \right), \cr
}
$$
\pagelabel{\Qfivesix}
where $\Lambda$ varies over all the subsets of $\{1,2,3,4\}$,
and $\C_1 + \left(b_{03},b_{04}\right)$.
However,
the latter prime ideal is not associated to $K(n,d)$ as
$K(n,d) : \bigcup_{i=1}^4 (b_{02}-b_{1i}b_{03})
= \C_1 + (b_{01}b_{03}^d - b_{04} b_{02}^d)$.
Also,
$Q_{1\emptyset}$ is not associated for the same reason.
On the other hand,
$Q_{1\Lambda}$ is associated to $K(n,d)$ when $\Lambda \not = \emptyset$
as it is minimal over it.
Thus

\prop
The set of associated primes of $K(n,d)$ is contained in
$\{Q_{1\Lambda'} \} \cup
\Ass \left( {R \over K(n,d) : b_{04}^d} \right)$,
as $\Lambda'$ varies over non-empty subsets of $\{1,2,3,4\}$.
\qed
\endb

To get at the remaining associated primes of $K(n,d)$,
one needs to calculate $K(n,d) : b_{04}^d$.
First let $\KK'$ (respectively $\EE'$)
be the ideal obtained from $\KK$ (respectively $\EE$)
by rewriting each $c_{1i}b_{01}^d$
as $c_{1i}b_{1i}^{d^2} b_{04}^d$.
Then both $\KK'$ and $\EE'$ are multiples of $b_{04}^d$.
As $c_{1i}(b_{01} - b_{1i}^d b_{04}) \in \FF$,
it follows that $\KK'+\EE' + \FF = \KEF$.
Thus
$$
\eqalignno{
\KEF : b_{04}^d &=
\KK'/b_{04}^d + \EE'/b_{04}^d + (\FF : b_{04}^d) \cr
&= \KK'/b_{04}^d
+ (c_{11}-b_{12}^{d^2} c_{12},c_{14}-c_{11},
b_{13}^{d^2} c_{13}-b_{12}^{d^2} c_{12}) \cr
&\hskip1em+
(c_{12}b_{11}-c_{13}b_{14},
c_{12}c_{2i}(b_{12}-b_{2i}b_{13}),
b_{01}b_{03}^d-b_{04}b_{02}^d) \cr
&\hskip1em+ c_{1i}(b_{02}-b_{1i}b_{03},
c_{1j}(b_{1i}-b_{1j}), b_{01}-b_{1i}^db_{04}). \cr
}
$$
By Fact~\factassocses,
the associated primes of $\KEF : b_{04}^d$
are in the union of the associated primes of the two ideals
obtained from $\KEF : b_{04}^d$
by respectively adding and coloning out $c_{12}$.
Note that $(\KEF: b_{04}^d)+ (c_{12})$ equals
$$
\eqalignno{
&= (c_{11}, c_{12},c_{14},
b_{01}b_{03}^d-b_{04}b_{02}^d) + c_{13}(b_{13}^{d^2},b_{14}, b_{02}-b_{13}b_{03},
b_{01}-b_{13}^db_{04}) \cr
&= (\C_1 + (b_{01}b_{03}^d-b_{04}b_{02}^d)) \cap
(c_{11}, c_{12},c_{14},b_{13}^{d^2}, b_{14},
b_{02}-b_{13}b_{03}, b_{01}-b_{13}^db_{04}), \cr
}
$$
whose associated prime ideals are
$$
\eqalignno{
Q_{2}&=\C_1 + (b_{01}b_{03}^d-b_{04}b_{02}^d), \cr
Q_{3}&= (c_{11}, c_{12},c_{14},b_{01},b_{02},b_{13},b_{14}). \cr
}
$$
Note that $Q_2$ is minimal over $K(n,d)$
and it is straightforward to verify that $Q_3$ is associated to $K(n,d)$.
Hence

\prop
\label{\propIrr}
$\displaystyle
\Ass \left( {R \over K(n,d) } \right) \subseteq
\{Q_{1\Lambda'}, Q_j| j = 2,3\} \cup
\Ass \left( {R \over K(n,d) : b_{04}^dc_{12}}\right)$,
where $\Lambda'$ varies over all the non-empty subsets of $\{1,2,3,4\}$.
\eqed
\endb

Let $\KK''$ be the ideal obtained from $\KK'/b_{04}^d$
by rewriting each $c_{11}$ and $c_{14}$
as $b_{12}^{d^2} c_{12}$.
Note that in the displayed $\KEF: b_{04}^d$ above,
the summand $\KK'/b_{04}^d$ may be replaced by $\KK''$,
and that
$\KK'' = (\KEpone) c_{12} b_{12}^{d^2}$.
Thus by Fact~\factcolony:
$$
\eqalignno{
&\hskip.5em\KEF: b_{04}^dc_{12}
= (\KEpone) b_{12}^{d^2} + (c_{11}-b_{12}^{d^2} c_{12},
c_{14}-c_{11},
c_{2i}(b_{12}-b_{2i}b_{13})) \cr
&\hskip1em+ (
b_{02}-b_{12}b_{03},
c_{1i}(b_{12}-b_{1i}),
b_{01}-b_{12}^db_{04}) \cr
&\hskip1em+ b_{12}^{d^2} (
b_{02}-b_{11}b_{03},
b_{02}-b_{14}b_{03},
b_{01}-b_{11}^db_{04},
b_{01}-b_{14}^db_{04}) \cr
&\hskip1em+ \bigl(
(
b_{01}b_{03}^d-b_{04}b_{02}^d,
b_{13}^{d^2} c_{13}-b_{12}^{d^2} c_{12},
c_{12}b_{11}-c_{13}b_{14}) + c_{13}(
b_{02}-b_{13}b_{03},
b_{01}-b_{13}^db_{04}) \bigr) : c_{12}. \cr
}
$$
The ideal in the last two rows,
before taking the colon with $c_{12}$,
decomposes as (by coloning and adding $c_{13}$):
$$
\eqalignno{
&\bigl(
b_{13}^{d^2} c_{13}-b_{12}^{d^2} c_{12},
c_{12}b_{11}-c_{13}b_{14},
b_{11}b_{13}^{d^2} -b_{14}b_{12}^{d^2},
b_{02}-b_{13}b_{03},
b_{01}-b_{13}^db_{04} \bigr) \cr
&\hskip1em \cap (
c_{13}, b_{12}^{d^2} c_{12},
c_{12}b_{11},
b_{01}b_{03}^d-b_{04}b_{02}^d), \cr
}
$$
which coloned with $c_{12}$ equals
$$
\eqalignno{
&\bigl(b_{13}^{d^2} c_{13}-b_{12}^{d^2} c_{12},
c_{12}b_{11}-c_{13}b_{14},
b_{11}b_{13}^{d^2} -b_{14}b_{12}^{d^2},
b_{02}-b_{13}b_{03}, b_{01}-b_{13}^db_{04} \bigr) \cr
&\hskip1em
\cap (c_{13}, b_{12}^{d^2}, b_{11},
b_{01}b_{03}^d-b_{04}b_{02}^d) \cr
=& (b_{13}^{d^2} c_{13}-b_{12}^{d^2} c_{12},
c_{12}b_{11}-c_{13}b_{14},
b_{11}b_{13}^{d^2} -b_{14}b_{12}^{d^2},b_{01}b_{03}^d-b_{04}b_{02}^d ) \cr
&\hskip2em
+ (b_{02}-b_{13}b_{03}, b_{01}-b_{13}^db_{04})
\cdot (c_{13}, b_{12}^{d^2}, b_{11}). \cr
}
$$
It follows that
$$
\eqalignno{
&\KEF: b_{04}^dc_{12}
= (\KEpone) b_{12}^{d^2}
+ (c_{11}-b_{12}^{d^2} c_{12},c_{14}-c_{11},
c_{2i}(b_{12}-b_{2i}b_{13})) \cr
&\hskip2em+ (b_{02}-b_{12}b_{03}, c_{1i}(b_{12}-b_{1i}),
b_{01}-b_{12}^db_{04}) \cr
&\hskip2em+ b_{12}^{d^2} (
b_{02}-b_{11}b_{03}, b_{02}-b_{14}b_{03},
b_{01}-b_{11}^db_{04}, b_{01}-b_{14}^db_{04}) \cr
&\hskip2em +(b_{13}^{d^2} c_{13}-b_{12}^{d^2} c_{12},
c_{12}b_{11}-c_{13}b_{14},
b_{11}b_{13}^{d^2} -b_{14}b_{12}^{d^2},b_{01}b_{03}^d-b_{04}b_{02}^d) \cr
&\hskip2em
+ (b_{02}-b_{13}b_{03}, b_{01}-b_{13}^db_{04})
\cdot (c_{13}, b_{12}^{d^2}, b_{11}) \cr
&= (\KEpone) b_{12}^{d^2}
+ (c_{11}-b_{12}^{d^2} c_{12},
c_{14}-c_{11}, c_{2i}(b_{12}-b_{2i}b_{13})) \cr
&\hskip2em+ (b_{02}-b_{12}b_{03},
c_{12}b_{12}^{d^2} (b_{12}-b_{11}),
c_{13}(b_{12}-b_{13}),
c_{12}b_{12}^{d^2} (b_{12}-b_{14}),
b_{01}-b_{12}^db_{04}) \cr
&\hskip2em+ b_{12}^{d^2} ((b_{12}-b_{1i})b_{03},
(b_{12}^d-b_{1i}^d)b_{04}, c_{13}-c_{12}) \cr
&\hskip2em +(c_{12}b_{11}-c_{13}b_{14},
b_{11}b_{13}^{d^2} -b_{14}b_{12}^{d^2})
+ b_{11}((b_{12}-b_{13})b_{03},
(b_{12}^d-b_{13}^d)b_{04}). \cr
}
$$
\pagelabel{\qKEFcolon}
When $n = 2$,
this ideal is much simpler than when $n > 2$.
We first analyze the cases $n \ge 3$.
By Fact~\factassocses,
the set of associated primes of $V$ is contained in
$\Ass \left( {R \over (\KEF : b_{04}^dc_{12}) + (b_{12}^{d^2})} \right)
\cup \Ass \left( {R \over \KEF : b_{04}^dc_{12} b_{12}^{d^2}} \right)$.
Note that
$$
\eqalignno{
&\hskip1.5em(\KEF : b_{04}^dc_{12}) + (b_{12}^{d^2}) = (c_{11}, c_{14},
b_{02}-b_{12}b_{03}, b_{01}-b_{12}^db_{04},
b_{12}^{d^2},
c_{13}(b_{12}-b_{13})) \cr
&\hskip2em+  (c_{12}b_{11}-c_{13}b_{14},c_{2i}(b_{12}-b_{2i}b_{13}))
+b_{11}(b_{13}^{d^2},(b_{12}-b_{13})b_{03},
(b_{12}^d-b_{13}^d)b_{04}). \cr
}
$$
\pagelabel{\qV}
To decompose this,
as Fact~\factfewvars\ applies here,
for notational simplicity it suffices to decompose
$$
V' =
(b_{12}^{d^2},
c_{13}(b_{12}-b_{13}),
c_{12}b_{11}-c_{13}b_{14},
c_{2i}(b_{12}-b_{2i}b_{13}))
+b_{11}(b_{13}^{d^2},(b_{12}-b_{13})b_{03},
(b_{12}^d-b_{13}^d)b_{04}).
$$
Then by first coloning and adding $c_{13}$, then repeating with $b_{03}b_{11}$:
$$
\eqalignno{
V' &= (b_{12}^{d^2}, c_{2i}b_{12}(1-b_{2i}),
b_{12}-b_{13}, c_{12}b_{11}-c_{13}b_{14}) \cr
&\hskip1em \cap \bigl((b_{12}^{d^2},
c_{2i}(b_{12}-b_{2i}b_{13}), c_{13})
+ b_{11}(
c_{12}, b_{13}^{d^2},(b_{12}-b_{13})b_{03},
(b_{12}^d-b_{13}^d)b_{04}) \bigr) \cr
&=(b_{12}^{d^2}, c_{2i}b_{12}(1-b_{2i}),
b_{12}-b_{13}, c_{12}b_{11}-c_{13}b_{14}) \cr
&\hskip1em \cap (b_{12}^{d^2}, c_{2i}b_{12}(1-b_{2i}),
c_{13}, c_{12}, b_{12}-b_{13}) \cr
&\hskip1em \cap \bigl((b_{12}^{d^2},
c_{2i}(b_{12}-b_{2i}b_{13}), c_{13})
+ b_{11}(c_{12}, b_{13}^{d^2}, b_{03},
(b_{12}^d-b_{13}^d)b_{04}) \bigr). \cr
}
$$
As the second component contains the first one,
the second one is redundant.
By decomposing the remaining two components of $V'$,
$V'$ further decomposes as:
$$
\eqalignno{
V'\kern-1ex&=(b_{12}^{d^2}, c_{2i}(1-b_{2i}),
b_{12}-b_{13}, c_{12}b_{11}-c_{13}b_{14}) \cr
&\hskip1em\cap (b_{12}, b_{13}, c_{12}b_{11}-c_{13}b_{14}) \cr
&\hskip1em \cap (b_{12}^{d^2},
c_{2i}(b_{12}-b_{2i}b_{13}), c_{13},b_{11}) \cr
&\hskip1em \cap \Bigl((b_{12}^{d^2},
c_{2i}(b_{12}-b_{2i}b_{13}),
c_{12}, c_{13}, b_{13}^{d^2}, b_{03},
(b_{12}^d-b_{13}^d)b_{04}) \Bigr). \cr
}
$$
Thus $V'$  is the intersection of the four ideals above.
Let the $i$th ideal be $V_i$.
The ideal $V_1$ decomposes into primary ideals as follows:
$$
\eqalignno{
V_1 &= \bigcap_{\Lambda} \Bigl((b_{12}^{d^2},
b_{12}-b_{13}, c_{12}b_{11}-c_{13}b_{14})
+ (c_{2i} | i \not \in \Lambda)
+ (1-b_{2i} | i \in \Lambda) \Bigr), \cr
}
$$
$V_2$ is a prime ideal,
$V_3$ decomposes as
$$
\eqalignno{
V_3 &=
\bigcap_{\Lambda} \Bigl((b_{12}^{d^2}, c_{13},b_{11}) +
(c_{2i} | i \not \in \Lambda) +
(b_{12}-b_{2i}b_{13} | i \in \Lambda) \Bigr) \cr
&=
\bigcap_{\Lambda} \Bigl((b_{12}^{d^2}, c_{13},b_{11})+ (c_{2i} | i \not \in \Lambda) +
(b_{12}-b_{2i}b_{13},b_{2i}^{d^2},
b_{2i}-b_{2i} | i,j \in \Lambda) \Bigr) \cr
&\hskip1em
\bigcap_{\Lambda \not = \emptyset} \Bigl((b_{12}^{d^2}, b_{13}^{d^2},
c_{13},b_{11}) + (c_{2i} | i \not \in \Lambda) +
(b_{12}-b_{2i}b_{13}| i \in \Lambda) \Bigr) \cr
&=
\bigcap_{\Lambda} \Bigl((b_{12}^{d^2}, c_{13},b_{11})+ (c_{2i} | i \not \in \Lambda) +
(b_{12}-b_{2i}b_{13},b_{2i}^{d^2},
b_{2i}-b_{2j} | i,j \in \Lambda) \Bigr) \cr
&\hskip1em
\bigcap_{\Lambda \not = \emptyset}
\Bigl((b_{12}^{d^2}, b_{13}^{d^2},c_{13},b_{11})+ (c_{2i} | i \not \in \Lambda) +
(b_{12}-b_{2i}b_{13}, b_{2i}-b_{2j} | i,j \in \Lambda) \Bigr) \cr
&\hskip1em \cap (b_{12}, b_{13}, c_{13},b_{11}), \cr
}
$$
which are all primary ideals.
The last primary factor of $V_3$ properly contains $V_2$
and is thus redundant in a primary decomposition of $V'$.
Finally,
$$
\eqalignno{
V_4 &=
\cap (b_{12}^{d^2}, c_{2i}(b_{12}-b_{2i}b_{13}),
c_{12}, c_{13}, b_{13}^{d^2}, b_{03},
b_{12}^d-b_{13}^d) \cr
&\hskip1em \cap (b_{12}^{d^2},
c_{2i}(b_{12}-b_{2i}b_{13}),
c_{12}, c_{13}, b_{13}^{d^2}, b_{03}, b_{04}) \cr
&= \bigcap_{\Lambda} \Bigl((b_{12}^{d^2},
c_{12}, c_{13}, b_{03}, b_{12}^d-b_{13}^d)+ (c_{2i} | i \not \in \Lambda)
+ (b_{12}-b_{2i}b_{13},
b_{2i}-b_{2j},1-b_{2i}^d | i,j \in \Lambda) \Bigr) \cr
&\hskip1em \bigcap_{\Lambda \not = \emptyset} \Bigl((b_{12}^d,
c_{12}, c_{13}, b_{13}^d, b_{03})+ (c_{2i} | i \not \in \Lambda)
+ (b_{12}-b_{2i}b_{13}, b_{2i}-b_{2j} | i,j \in \Lambda)
\Bigr) \cr
&\hskip1em \cap\Bigl(b_{12},b_{13},c_{12}, c_{13}, b_{03}\Bigr) \cr
&\hskip1em \bigcap_{\Lambda} \Bigl((b_{12}^{d^2},
c_{12}, c_{13}, b_{13}^{d^2}, b_{03}, b_{04})+ (c_{2i} | i \not \in \Lambda) +
(b_{12}-b_{2i}b_{13},b_{2i}-b_{2j} | i,j \in \Lambda)
\Bigr) \cr
&\hskip1em \cap (b_{12},b_{13}, c_{12}, c_{13}, b_{03}, b_{04}). \cr
}
$$
The ideals in the last and the third to the last rows above contain $V_2$,
and are thus redundant in a primary decomposition of $V$.
Now the (possibly redundant) associated primes of $V'$ can be easily read off,
and by adding $(c_{11},c_{14},b_{02}-b_{12}b_{03},
b_{01}-b_{12}^db_{04})$
one gets the associated primes of
$(\KEF: b_{04}^dc_{12}) + (b_{12}^{d^2})$
(see page~\qV):
$$
\eqalignno{
Q_{4\Lambda}&=
(c_{11}, c_{14}, b_{01}, b_{02}, b_{12},b_{13}, c_{12}b_{11}-c_{13}b_{14})
+ (c_{2i} | i \not \in \Lambda) + (1-b_{2i} | i \in \Lambda), \cr
Q_{5}&=
(c_{11}, c_{14}, b_{01}, b_{02}, b_{12}, b_{13}, c_{12}b_{11}-c_{13}b_{14}), \cr
Q_{6\Lambda}&=
(c_{11}, c_{13}, c_{14}, b_{01}, b_{02}, b_{11}, b_{12})
+ (c_{2i} | i \not \in \Lambda) + (b_{2i} | i \in \Lambda), \cr
Q_{7\Lambda'}&=
(c_{11}, c_{13}, c_{14}, b_{01}, b_{02}, b_{11}, b_{12},b_{13})
+(c_{2i} | i \not \in \Lambda') +(b_{2i}-b_{2j} | i,j \in \Lambda'), \cr
Q_{8\Lambda\alpha}&=
\C_1 + (b_{01}, b_{02}, b_{03}, b_{12},b_{13})
+ (c_{2i} | i \not \in \Lambda) + (b_{2i}-\alpha | i \in \Lambda), \alpha^d = 1, \cr
Q_{9\Lambda'}&=
\C_1 + (b_{01}, b_{02}, b_{03}, b_{12},b_{13})
+ (c_{2i} | i \not \in \Lambda')
+ (b_{2i}-b_{2j} | i,j \in \Lambda'), \cr
Q_{10\Lambda}&=\C_1 +
(b_{01}, b_{02}, b_{03}, b_{04}, b_{12}, b_{13})
+ (c_{2i} | i \not \in \Lambda) +
(b_{2i}-b_{2j} | i,j \in \Lambda), \cr
}
$$
where $\Lambda$ varies over all subsets of $\{1,2,3,4\}$
and $\Lambda'$ varies over all the non-empty subsets of $\{1,2,3,4\}$.
By Proposition~\propIrr\ it follows that

\prop
\label{\propIrrr}
Whenever $n \ge 3$,
$$
\eqalignno{
\Ass &\left( {R \over K(n,d) } \right) \subseteq
\{Q_{i\Lambda}, Q_j, Q_{k\Lambda'}|
i = 4, 6, 10; j =  2,3,5; k = 1, 7, 9 \} \cr
& \cup \{Q_{8\Lambda\alpha}\} \cup
\Ass \left( {R \over K(n,d) :
b_{04}^dc_{12}b_{12}^{d^2}}\right), \cr
}
$$
where $\alpha$ varies over the $d$th roots of unity
and $\Lambda$ over all the subsets of $\{1,2,3,4\}$
and $\Lambda'$ varies over all the non-empty subsets of $\{1,2,3,4\}$.\qed
\endb

Next,
\pagelabel{\qcalco}
$$
\eqalignno{
&\KEF: b_{04}^dc_{12}b_{12}^{d^2}= \KEpone
+ (c_{11}-b_{12}^{d^2}c_{12},
c_{14}-c_{11}, b_{02}-b_{12}b_{03},
b_{01}-b_{12}^db_{04}) \cr
&\hskip2em+ (
c_{12}(b_{12}-b_{11}),
c_{12}(b_{12}-b_{14}),
(b_{12}-b_{1i})b_{03}),
(b_{12}^d-b_{1i}^d)b_{04},
c_{13}-c_{12}) \cr
&\hskip2em+ \Bigl( (
c_{2i}(b_{12}-b_{2i}b_{13}),
c_{13}(b_{12}-b_{13}),
c_{12}b_{11}-c_{13}b_{14},
b_{11}b_{13}^{d^2}-b_{14}b_{12}^{d^2})
\cr
&\hskip3em + b_{11}(
(b_{12}-b_{13})b_{03},
(b_{12}^d-b_{13}^d)b_{04})
\Bigr) : b_{12}^{d^2}. \cr
}
$$
The ideal $\hat V$ in the last two rows,
before coloning with $b_{12}^{d^2}$,
(partially) decomposes as (first coloning and adding $c_{13}$):
$$
\eqalignno{
\hat V&= (c_{2i}b_{12}(1-b_{2i}), b_{12}-b_{13},
c_{12}b_{11}-c_{13}b_{14},
b_{12}^{d^2}(b_{11}-b_{14})) \cr
&\hskip1em \cap \bigl(c_{2i}(b_{12}-b_{2i}b_{13}), c_{13},
c_{12}b_{11},
b_{11}b_{13}^{d^2}-b_{14}b_{12}^{d^2},(b_{12}-b_{13})b_{03}b_{11},
(b_{12}^d-b_{13}^d)b_{04}b_{11}) \bigr). \cr
}
$$
The second component decomposes further into:
$$
\eqalignno{
&\bigcap_{\Lambda} \Bigl((c_{13}, c_{12}b_{11},
b_{11}b_{13}^{d^2}-b_{14}b_{12}^{d^2},
(b_{12}-b_{13})b_{03}b_{11},(b_{12}^d-b_{13}^d)b_{04}b_{11}) \cr
&\hskip2em
+ (c_{2i} | i \not \in \Lambda)
+ (b_{12}-b_{2i}b_{13} | i \in \Lambda) \Bigr). \cr
}
$$
For $\Lambda \not = \emptyset$,
the corresponding component of $\hat V$ simplifies to
$$
\eqalignno{
&(c_{13}, c_{12}b_{11}) + (c_{2i} | i \not \in \Lambda) \cr
&\hskip1em
+ (b_{12}-b_{2i}b_{13},
b_{13}^{d^2}(b_{11}-b_{14}b_{2i}^{d^2}),b_{13}(b_{2i}-1)b_{03}b_{11},
b_{13}^d(b_{2i}^d-1)b_{04}b_{11}
| i \in \Lambda \Bigr), \cr
}
$$
and for $\Lambda = \emptyset$ it equals to and (partially) decomposes as
$$
\eqalignno{
&\C_2 + (c_{13}, c_{12}b_{11},
b_{11}b_{13}^{d^2}-b_{14}b_{12}^{d^2},
(b_{12}-b_{13})b_{03}b_{11},(b_{12}^d-b_{13}^d)b_{04}b_{11}) \cr
&= \bigl( \C_2 + (c_{13}, b_{11}, b_{14}b_{12}^{d^2})
\bigr) \cr
&\hskip1em \cap
\bigl(\C_2 + (c_{13}, c_{12},
b_{11}b_{13}^{d^2}-b_{14}b_{12}^{d^2},
(b_{12}-b_{13})b_{03}b_{11},(b_{12}^d-b_{13}^d)b_{04}b_{11})\bigr) \cr
&= \bigl( \C_2 + (c_{13}, b_{11}, b_{14}b_{12}^{d^2})
\bigr) \cr
&\hskip1em \cap
\bigl(\C_2 + (c_{13}, c_{12},
(b_{11}-b_{14})b_{12}^{d^2},
b_{12}-b_{13}) \bigr) \cr
&\hskip1em \cap
\bigl(\C_2 + (c_{13}, c_{12},b_{11}b_{13}^{d^2}-b_{14}b_{12}^{d^2},
b_{03}b_{11},
(b_{12}^d-b_{13}^d)b_{04}b_{11}). \cr
}
$$
The last component in the last display above decomposes as:
$$
\eqalignno{
&\bigl(\C_2 + (c_{13}, c_{12},
(b_{11}-b_{14})b_{12}^{d^2},
b_{03}, b_{12}^d-b_{13}^d\bigr) \cr
&\hskip1em \cap \bigl(\C_2 + (c_{13}, c_{12},
b_{11}b_{13}^{d^2}-b_{14}b_{12}^{d^2},
b_{03}, b_{04})\bigr) \cr
&\hskip1em \cap \bigl(\C_2 + (c_{13}, c_{12},
b_{14}b_{12}^{d^2}, b_{11})\bigr). \cr
}
$$
Thus from the combined decomposition above of $\hat V$,
$\hat V : b_{12}^{d^2}$ equals
$$
\eqalignno{
\hat V : b_{12}^{d^2}
&= (c_{2i}(1-b_{2i}), b_{12}-b_{13},
c_{12}b_{11}-c_{13}b_{14}, b_{11}-b_{14}) \cr
&\hskip1em\cap \bigcap_{\Lambda \not = \emptyset} \bigl(
(c_{13}, c_{12}b_{14}) + (c_{2i} | i \not \in \Lambda)
+ (b_{12}-b_{2i}b_{13}, b_{2i}-b_{2j}
| i, j\in \Lambda) \cr
&\hskip2em + (b_{11}-b_{14}b_{2i}^{d^2},
(b_{2i}-1)b_{03}b_{14},
(b_{2i}^d-1)b_{04}b_{14}
| i \in \Lambda) \bigr) \cr
&\hskip1em \cap
\bigl( \C_2 + (c_{13}, b_{11}, b_{14})
\bigr) \cr
&\hskip1em \cap
\bigl(\C_2 + (c_{13}, c_{12}, b_{11}-b_{14},
b_{12}-b_{13}) \bigr) \cr
&\hskip1em \cap \bigl(\C_2 + (c_{13}, c_{12}, b_{11}-b_{14},
b_{03}, b_{12}^d-b_{13}^d\bigr) \cr
&\hskip1em \cap \bigl(\C_2 + (c_{13}, c_{12},
b_{11}b_{13}^{d^2}-b_{14}b_{12}^{d^2},
b_{03}, b_{04})\bigr) \cr
&\hskip1em \cap \bigl(\C_2 + (c_{13}, c_{12},
b_{14}, b_{11})\bigr). \cr
}
$$
Next we compute this intersection.
The intersection of the last two ideals is
$$
\C_2 + (c_{13}, c_{12},
b_{11}b_{13}^{d^2}-b_{14}b_{12}^{d^2},
b_{03}b_{11}, b_{04}b_{11},b_{03}b_{14}, b_{04}b_{14}),
$$
and the intersection of the two immediately preceding is
$$
\C_2 + (c_{13}, c_{12}, b_{11}-b_{14},
b_{03}(b_{12}-b_{13}),
b_{12}^d-b_{13}^d\bigr).
$$
Thus the intersection of the last four components in $\hat V:b_{12}^{d^2}$ is
$$
\eqalignno{
\C_2 &+ (c_{13}, c_{12},
b_{11}b_{13}^{d^2}-b_{14}b_{12}^{d^2},
b_{03}(b_{11}-b_{14}), b_{04}(b_{11}-b_{14}) \cr
&+(b_{03}b_{11}(b_{12}-b_{13}),
b_{04}b_{14}(b_{12}^d-b_{13}^d)), \cr
}
$$
and then the intersection of the last five components in $\hat V:b_{12}^{d^2}$ is
$$
\eqalignno{
\C_2 &+ (c_{13},
b_{11}b_{13}^{d^2}-b_{14}b_{12}^{d^2},
b_{03}(b_{11}-b_{14}), b_{04}(b_{11}-b_{14}) \cr
&+(b_{03}b_{11}(b_{12}-b_{13}),
b_{04}b_{14}(b_{12}^d-b_{13}^d))+
c_{12}(b_{11},b_{14}). \cr
}
$$
Now this ideal,
together with the second through the sixteenth components of $\hat V:b_{12}^{d^2}$
above are all of the form
$$
\eqalignno{
(c_{13}, &b_{11}b_{13}^{d^2}-b_{14}b_{12}^{d^2},
b_{03}(b_{11}-b_{14}), b_{04}(b_{11}-b_{14}),b_{03}b_{11}(b_{12}-b_{13}),
b_{04}b_{14}(b_{12}^d-b_{13}^d)) \cr
&\hskip2em+
c_{12}(b_{11},b_{14})+ (c_{2i} | i \not \in \Lambda) \cr
&\hskip2em
+ (b_{12}-b_{2i}b_{13}, b_{2i}-b_{2j}
(b_{11}-b_{14}b_{2i}^{d^2},
(b_{2i}-1)b_{03}b_{14},
(b_{2i}^d-1)b_{04}b_{14}
| i, j\in \Lambda) \bigr), \cr
}
$$
as $\Lambda$ varies over all the sixteen subsets of $\{1,2,3,4\}$.
The intersection of all these sixteen ideals is
$$
\eqalignno{
&\hskip1em(c_{13}, b_{11}b_{13}^{d^2}-b_{14}b_{12}^{d^2},
b_{03}(b_{11}-b_{14}), b_{04}(b_{11}-b_{14}),b_{03}b_{11}(b_{12}-b_{13}),
b_{04}b_{14}(b_{12}^d-b_{13}^d)) \cr
&+
c_{12}(b_{11},b_{14})+ c_{2i} (b_{12}-b_{2i}b_{13}, c_{2j}(b_{2i}-b_{2j}),
b_{11}-b_{14}b_{2i}^{d^2}),(b_{2i}-1)b_{03}b_{14},
(b_{2i}^d-1)b_{04}b_{14})). \cr
}
$$
Thus finally $\hat V : b_{12}^{d^2}$ equals
(intersection of above with the first factor of $\hat V :b_{12}^{d^2}$):
$$
\eqalignno{
(&b_{11}b_{13}^{d^2}-b_{14}b_{12}^{d^2},
b_{03}(b_{11}-b_{14}), b_{04}(b_{11}-b_{14}),b_{03}b_{11}(b_{12}-b_{13}),
b_{04}b_{14}(b_{12}^d-b_{13}^d)) \cr
&+ c_{2i} (b_{12}-b_{2i}b_{13}, c_{2j}(b_{2i}-b_{2j}),
b_{11}-b_{14}b_{2i}^{d^2},(b_{2i}-1)b_{03}b_{14},
(b_{2i}^d-1)b_{04}b_{14})) \cr
&+
(c_{12}(b_{11}-b_{14},c_{12}b_{11}-c_{13}b_{14})
+c_{13}(c_{2i}(1-b_{2i}), b_{12}-b_{13},b_{11}-b_{14}). \cr
}
$$
Thus $\KEF: b_{04}^dc_{12}b_{12}^{d^2}$ equals
(refer to page~\qcalco):
$$
\eqalignno{
&\KEpone
+ (c_{11}-b_{12}^{d^2}c_{12},
c_{14}-c_{11}, b_{02}-b_{12}b_{03},
b_{01}-b_{12}^db_{04} ) \cr
&\hskip2em+ (c_{12}(b_{12}-b_{11}),
c_{12}(b_{12}-b_{14}),(b_{12}-b_{1i})b_{03}, (b_{12}^d-b_{1i}^d)b_{04},
c_{13}-c_{12} ) \cr
&\hskip2em+(b_{11}b_{13}^{d^2}-b_{14}b_{12}^{d^2},
b_{03}(b_{11}-b_{14}), b_{04}(b_{11}-b_{14}) \cr
&\hskip2em+(b_{03}b_{11}(b_{12}-b_{13}),
b_{04}b_{14}(b_{12}^d-b_{13}^d)) \cr
&\hskip2em+ c_{2i} (b_{12}-b_{2i}b_{13}, c_{2j}(b_{2i}-b_{2j}),
b_{11}-b_{14}b_{2i}^{d^2},(b_{2i}-1)b_{03}b_{14},
(b_{2i}^d-1)b_{04}b_{14})) \cr
&\hskip2em+(c_{12}(b_{11}-b_{14}),c_{12}b_{11}-c_{13}b_{14})
+c_{13}(c_{2i}(1-b_{2i}), b_{12}-b_{13},b_{11}-b_{14}) \cr
&= \KEpone
+ (c_{11}-b_{12}^{d^2}c_{12},
c_{14}-c_{11}, c_{13}-c_{12}, b_{02}-b_{12}b_{03},b_{01}-b_{12}^db_{04} ) \cr
&\hskip2em+ (
c_{12}(b_{12}-b_{1i}),
(b_{12}-b_{1i})b_{03},
(b_{12}^d-b_{1i}^d)b_{04},
b_{11}b_{13}^{d^2}-b_{14}b_{12}^{d^2},
(b_{11}-b_{14})b_{04}) \cr
&\hskip2em+ c_{2i}(
b_{12}-b_{2i}b_{13}, c_{2j}(b_{2i}-b_{2j}),
b_{11}-b_{2i}^{d^2}b_{14},
c_{12}(1-b_{2i}) ). \cr
}
$$
\pagelabel{\qKEFVcolon}
By Fact~\factassocses,
any associated prime of $\KEF : b_{04}^dc_{12}b_{12}^{d^2}$
is associated either to
$(\KEF : b_{04}^dc_{12}b_{12}^{d^2})+ (c_{12})$
or to $\KEF : b_{04}^dc_{12}^2b_{12}^{d^2}$.
These ideals are as follows:
$$
\eqalignno{
&A= \KEF : b_{04}^dc_{12}^2 b_{12}^{d^2}= \KEpone \cr
&
+ (c_{11}-b_{12}^{d^2}c_{12},
c_{14}-c_{11}, c_{13}-c_{12},b_{02}-b_{12}b_{03},
b_{12}-b_{1i},
b_{01}-b_{12}^db_{04},
c_{2i}(1-b_{2i}) ), \cr
}
$$
and
$$
\eqalignno{
&(\KEF : b_{04}^dc_{12}b_{12}^{d^2})+ (c_{12})
= \KEpone + \C_1 \cr
&\hskip2em+ (b_{02}-b_{12}b_{03},
b_{01}-b_{12}^db_{04},
(b_{12}-b_{1i})b_{03},
(b_{12}^d-b_{1i}^d)b_{04} ) \cr
&\hskip2em+ c_{2i}(
b_{12}-b_{2i}b_{13}, c_{2j}(b_{2i}-b_{2j}),
b_{11}-b_{2i}^{d^2}b_{14} )+ (
b_{11}b_{13}^{d^2}-b_{14}b_{12}^{d^2},
(b_{11}-b_{14})b_{04} ). \cr
}
$$
This last ideal decomposes (coloning and adding $b_{03}$):
$$
\eqalignno{
&= \Bigl( \KEpone + \C_1+ (b_{02}-b_{12}b_{03},
b_{01}-b_{12}^db_{04},
b_{12}-b_{1i} )+ c_{2i}(
b_{12}(1-b_{2i}), c_{2j}(b_{2i}-b_{2j}) ) \Bigr) \cr
&\hskip1em \cap \Bigl( \KEpone + \C_1
+ (b_{01}-b_{12}^db_{04},b_{02}, b_{03},
(b_{12}^d-b_{1i}^d)b_{04} ) \cr
&\hskip2em+ c_{2i}(
b_{12}-b_{2i}b_{13}, c_{2j}(b_{2i}-b_{2j}),
b_{11}-b_{2i}^{d^2}b_{14} )+ (
b_{11}b_{13}^{d^2}-b_{14}b_{12}^{d^2},
(b_{11}-b_{14})b_{04} ) \Bigr). \cr
}
$$
Let $B$ denote the first of the two components above.
The second component further decomposes
(coloning and adding $b_{04}$):
$$
\eqalignno{
&\Bigl( \KEpone + \C_1
+ (b_{01}-b_{12}^db_{04},b_{02}, b_{03},
b_{12}^d-b_{1i}^d,
b_{11}-b_{14}) \cr
&\hskip2em+ c_{2i}(
b_{12}-b_{2i}b_{13}, c_{2j}(b_{2i}-b_{2j}),
b_{11}-b_{2i}^{d^2}b_{14} ) \Bigr) \cr
&\hskip1em \cap \Bigl( \KEpone + \C_1
+ (b_{01}, b_{02}, b_{03}, b_{04},
b_{11}b_{13}^{d^2}-b_{14}b_{12}^{d^2} ) \cr
&\hskip2em+ c_{2i}(
b_{12}-b_{2i}b_{13}, c_{2j}(b_{2i}-b_{2j}),
b_{11}-b_{2i}^{d^2}b_{14} ) \Bigr). \cr
}
$$
Let these two ideals be $C$ and $D$,
in that order.
By Fact~\factassocses\ and Proposition~\propIrrr:

\prop
\label{\propIrrrr}
When $n \ge 3$,
$$
\eqalignno{
\Ass &\left( {R \over K(n,d) } \right) \subseteq
\{Q_{ir\Lambda}, Q_{jr},Q_{kr\Lambda'}, Q_{8\Lambda\alpha}
| i = 5, 6, 10; j =  2,3,5; k = 1, 7, 9 \} \cr
&\cup
\Ass \left( {R \over A} \right)
\cup \Ass \left( {R \over B} \right)
\cup \Ass \left( {R \over C} \right)
\cup \Ass \left( {R \over D} \right), \cr
}
$$
where $\alpha$ varies over all the $d$th roots of unity in $K$
and $\Lambda, \Lambda'$ over all the subsets of $\{1,2,3,4\}$
with $|\Lambda'| > 0$.
\qed
\endb

But $D = K(n-1,d^2) + \C_1 + (b_{01}, b_{02}, b_{03}, b_{04})$,
so the induction on $n$ gives all of its associated primes.
It remains to compute the associated primes of the other three ideals,
namely of $A$, $B$ and $C$.
For each one of these ideals $U$ ($U$ varying over $A, B$ and $C$),
by Fact~\factassocses,
$$
\Ass \left( {R \over U} \right) \subseteq
\Ass \left( {R \over U + (b_{11}^{d^2})} \right)
\cup \Ass \left( {R \over U : b_{11}^{d^2}} \right).
$$
Modulo the other generators of $U$,
$\KEpone$ is contained in $(b_{11}^{d^2})$,
so that the three $U + (b_{11}^{d^2})$ are as follows:
$$
\eqalignno{
A &+ (b_{11}^{d^2})=
(b_{11}^{d^2}, c_{11},
c_{14}-c_{11}, c_{13}-c_{12},
b_{02}-b_{12}b_{03},
b_{12}-b_{1i},
b_{01}-b_{12}^db_{04},
c_{2i}(1-b_{2i}) ), \cr
B &+ (b_{11}^{d^2})=
\C_1+ (b_{11}^{d^2}, b_{02}-b_{12}b_{03},
b_{01}-b_{12}^db_{04},
b_{12}-b_{1i} )+ c_{2i}(
b_{12}(1-b_{2i}), c_{2j}(b_{2i}-b_{2j}) ), \cr
C &+ (b_{11}^{d^2})=
\C_1
+ (b_{11}^{d^2},
b_{01}-b_{12}^db_{04},b_{02}, b_{03},
b_{12}^d-b_{1i}^d,
b_{11}-b_{14}) \cr
&\hskip2em+ c_{2i}(
b_{12}-b_{2i}b_{13}, c_{2j}(b_{2i}-b_{2j}),
b_{11}(1-b_{2i}^{d^2})). \cr
}
$$
The contribution of $A$ and $B$ to the possible associated primes of $K(n,d)$
are as follows:
$$
\eqalignno{
Q_{11\Lambda} &= (c_{11},c_{13}-c_{12}, c_{14},
b_{01},b_{02},b_{11},b_{12},b_{13},b_{14})
+ (c_{2i} | i \not \in \Lambda)
+ (1-b_{2i} | i \in \Lambda)  \hbox{(from $A$)}, \cr
Q_{12\Lambda} &= \C_1 +
(b_{01},b_{02},b_{11},b_{12},b_{13},b_{14})
+ (c_{2i} | i \not \in \Lambda)+ (1-b_{2i} | i\in \Lambda)
\hbox{\ \ (from $B$)}, \cr
Q_{13\Lambda} &= \C_1 +
(b_{01},b_{02},b_{11},b_{12},b_{13},b_{14})
+ (c_{2i} | i \not \in \Lambda)
+ (b_{2i}-b_{2j} | i, j \in \Lambda) \hbox{\ \ (from $B$)}. \cr
}
$$
The associated primes obtained from $C$ are not so easily read off,
so we need to first decompose $C + (b_{11}^{d^2})$.
By Fact~\factfewvars\ it suffices to decompose
$$
\eqalignno{
&\hskip1em(b_{11}^{d^2},
b_{12}^d-b_{11}^d,
b_{12}^d-b_{13}^d)+ c_{2i}(
b_{12}-b_{2i}b_{13}, c_{2j}(b_{2i}-b_{2j}),
b_{11}(1-b_{2i}^{d^2})) \cr
&= \bigcap_{\Lambda} \bigl((b_{11}^{d^2},
b_{12}^d-b_{11}^d,
b_{12}^d-b_{13}^d)+ (c_{2i} | i \not \in \Lambda)
+(b_{12}-b_{2i}b_{13}, b_{2i}-b_{2j},
b_{11}(1-b_{2i}^{d^2}) | i,j \in \Lambda) \bigr) \cr
&= \bigl(\C_2  + (b_{11}^{d^2},
b_{12}^d-b_{11}^d,
b_{12}^d-b_{13}^d)\bigr) \cr
&\hskip1em\cap \bigcap_{\Lambda\not = \emptyset} \bigl((b_{11}^{d^2},
b_{13}^d-b_{11}^d) + (c_{2i} | i \not \in \Lambda)\cr
&\hskip2em
+(b_{12}-b_{2i}b_{13}, b_{2i}-b_{2j},
b_{11}(1-b_{2i}^{d^2}),
b_{13}^d(b_{2i}^d-1)| i,j \in \Lambda) \bigr) \cr
&= \bigl(\C_2  + (b_{11}^{d^2},
b_{12}^d-b_{11}^d,
b_{12}^d-b_{13}^d)\bigr) \cr
&\hskip1em\bigcap_{\Lambda\not = \emptyset} \bigl((b_{11},
b_{13}^d) + (c_{2i} | i \not \in \Lambda)
+(b_{12}-b_{2i}b_{13}, b_{2i}-b_{2j}| i,j \in \Lambda) \bigr) \cr
&\hskip1em \bigcap_{\Lambda\not = \emptyset} \bigl((b_{11}^d,
b_{13}^d) + (c_{2i} | i \not \in \Lambda)
+(b_{12}-b_{2i}b_{13}, b_{2i}-b_{2j},
(1-b_{2i}^{d^2})/(b_{2i}^d-1)| i,j \in \Lambda) \bigr) \cr
&\hskip1em\cap \bigcap_{\Lambda\not = \emptyset} \bigl((b_{11}^{d^2},
b_{13}^d-b_{11}^d) + (c_{2i} | i \not \in \Lambda)
+(b_{12}-b_{2i}b_{13}, b_{2i}-b_{2j},
1-b_{2i}^d| i,j \in \Lambda) \bigr), \cr
}
$$
from which one can easily read off the associated primes.
Then by adding $\C_1+ (b_{01}-b_{12}^db_{04},b_{02}, b_{03},b_{11}-b_{14})$
to these primes,
the contribution of $C + (b_1^{d^2})$ to the associated primes of $K(n,d)$:
$$
\eqalignno{
Q_{14\Lambda} &= \C_1 +
(b_{01},b_{02},b_{03},b_{11},b_{12},b_{13},b_{14})
+ (c_{2i} | i \not \in \Lambda)
+ (b_{2i}- b_{2j} | i,j \in \Lambda), \cr
Q_{15\Lambda'\alpha} &= \C_1 +
(b_{01},b_{02},b_{03},b_{11},b_{12},b_{13},b_{14}) \cr
&\hskip2em
+ (c_{2i} | i \not \in \Lambda')
+ (b_{2i}- \alpha | i \in \Lambda'), \alpha^{d^2} = 1, \alpha^d \not= 1, \cr
Q_{16\Lambda'\alpha} &= \C_1 +
(b_{01},b_{02},b_{03},b_{11},b_{12},b_{13},b_{14})
+ (c_{2i} | i \not \in \Lambda')
+ (b_{2i}- \alpha | i \in \Lambda'), \alpha^d = 1, \cr
}
$$
where $\Lambda$ varies over all and $\Lambda'$ over all non-empty
subsets of $\{1,2,3,4\}$.
This establishes

\prop
\label{\propIrrrrr}
For $n \ge 3$,
$$
\eqalignno{
\Ass &\left( {R \over K(n,d)} \right) \subseteq
\{Q_{i\Lambda}, Q_j, Q_{k\Lambda'}|
i = 4, 6, 10,11,12,13,14; j = 2,3,5; k = 1, 7, 9\} \cr
&\cup
\{Q_{8\Lambda\alpha},Q_{15\Lambda'\alpha'},Q_{16\Lambda'\alpha}\}
\cup
\Ass \left( {R \over A:b_{11}^{d^2}} \right)
\cup \Ass \left( {R \over B:b_{11}^{d^2}} \right)
\cup \Ass \left( {R \over C:b_{11}^{d^2}} \right)
\cr
&\cup \Ass \left( {R \over K(n-1,d^2) + \C_1 + (b_{01}, b_{02}, b_{03}, b_{04})} \right), \cr
}
$$
where $\alpha$ varies over all the $d$th roots of unity in $K$,
$\alpha'$ varies over all the $d^2$ roots of unity in $K$
which are not $d$th roots of unity,
$\Lambda$ varies over all the subsets of $\{1,2,3,4\}$,
and $\Lambda'$ varies over all the non-empty subsets of $\{1,2,3,4\}$.
\qed
\endb

\vskip 4ex
Finally we analyze the three $U : b_{11}^{d^2}$
(where $U = A, B, C$, and $n \ge 3$).
For notational purposes we also define the following ideals in $R$:
$$
\eqalignno{
\D_r &= \left(c_{r4}-c_{r1},c_{r3}-c_{r2},c_{r2}-c_{r1} \right), r = 1, \ldots, n-1, \cr
\D_n &= (0), \cr
B_{0} &= B_{1} = (0), \cr
B_{r} &= \left(1-b_{2i}, 1-b_{3i}, \ldots, 1-b_{ri} |
i = 1, \ldots, 4 \right), r = 2, \ldots, n-1. \cr
B_{kr} &= \left(1-b_{ki}, 1-b_{k+1,i}, \ldots, 1-b_{ri} |
i = 1, \ldots, 4 \right), r = 2, \ldots, n-1. \cr
}
$$
Then the computation of $U : b_{11}^{d^2}$ is is not difficult,
as modulo $b_{11}^d-b_{14}^d \in U$,
$\KEpone$ equals $b_{11}^{d^2}\hat{\KK}$ for the following ideal $\hat{\KK}$:
$$
\eqalignno{
\hat{\KK} &=\bigl((g_{rj}| r \ge 4)
\hbox{\ with $c_{11}= 1$} \bigr)+ \D_2 +
(c_{22}(b_{21}-b_{24}),c_{22}c_{3i}(b_{12}-b_{2i}b_{13})). \cr
}
$$
In particular,
$\hat{\KK}$ contains $\D_2$.
Thus $A : b_{11}^{d^2}$ equals
$$
\hat{\KK}
+ (c_{11}-b_{12}^{d^2}c_{12},
c_{14}-c_{11}, c_{13}-c_{12},b_{02}-b_{12}b_{03},
b_{12}-b_{1i},
b_{01}-b_{12}^db_{04},
c_{22}(1-b_{2i}) )
$$
which decomposes as follows:
$$
\eqalignno{
&\bigcap_{t= 2} ^n \bigl(\D_2 + \cdots + \D_{t-1} + \C_t
+ B_{2,t-1}+ (c_{11}-b_{12}^{d^2}c_{12}) \cr
&\hskip2em
+ (c_{14}-c_{11}, c_{13}-c_{12},b_{02}-b_{12}b_{03},
b_{12}-b_{1i},b_{01}-b_{12}^db_{04} ) \bigr). \cr
}
$$
Thus $A : b_{11}^{d^2}$ contributes the following
possible associated prime ideals of $K(n,d)$,
for $t = 2 \ldots, n$:
$$
\eqalignno{
Q_{17t}&= \D_2 + \cdots + \D_{t-1} + \C_t
+ B_{2t-1} \cr
&\hskip2em+ (c_{11}-b_{12}^{d^2}c_{12},c_{14}-c_{11}, c_{13}-c_{12},b_{02}-b_{12}b_{03},
b_{12}-b_{1i},b_{01}-b_{12}^db_{04} ). \cr
}
$$
Similarly,
$B: b_{11}^{d^2}$ equals
$$
\hat{\KK}+ \C_1+
\left(b_{02}-b_{12}b_{03},
b_{01}-b_{12}^db_{04},
b_{12}-b_{1i},c_{22}(1-b_{2i}) \right),
$$
which contributes the following
possible associated prime ideals of $K(n,d)$ for $t = 2 \ldots, n$:
$$
\eqalignno{
Q_{18t}&= \C_1+ D_2 + \cdots + \D_{t-1} + \C_t+ B_{2t-1}+ \left( b_{02}-b_{12}b_{03},
b_{12}-b_{1i},
b_{01}-b_{12}^db_{04} \right). \cr
}
$$
And lastly,
$$
\eqalignno{
&C: b_{11}^{d^2} = \hat{\KK}+ \C_1
+ \left(b_{01}-b_{12}^db_{04},b_{02}, b_{03},
b_{11}-b_{14}\right) \cr
&\hskip2em+
\left(b_{12}^d-b_{1i}^d\right)
+ c_{22}\left(
b_{12}-b_{2i}b_{13}, b_{2i}-b_{2j},
1-b_{2i}^d) \right) \cr
&= \Bigl(\C_1+ \C_2
+ \left(b_{01}-b_{12}^db_{04},b_{02}, b_{03},
b_{11}-b_{14}, b_{12}^d-b_{1i}^d\right)\Bigr) \cr
&\hskip1em \bigcap_{t = 3}^n \Bigl(\C_1 + \D_2 + \cdots + \D_{t-1}
+ \C_t + B_{3t-1} \cr
&\hskip2em
+ \left(b_{01}-b_{12}^db_{04},b_{02}, b_{03},
b_{11}-b_{14},b_{12}^d-b_{1i}^d,
b_{12}-b_{2i}b_{13}, b_{2i}-b_{2j},
1-b_{2i}^d \right) \Bigr), \cr
}
$$
which contributes the following
possible associated prime ideals of $K(n,d)$ for $t = 3 \ldots, n$,
and with $\alpha^d = \beta^d = 1$:
$$
\eqalignno{
Q_{19,2\alpha\beta} &= \C_1+ \C_2
+ \left(b_{01}-b_{12}^db_{04},b_{02}, b_{03},
b_{12}-\alpha b_{13},b_{11}-\beta b_{13}, b_{11}-b_{14}
\right), \cr
Q_{19t\alpha\beta} &= \C_1+ \D_2 + \cdots + \D_{t-1}
+ \C_t + B_{3t-1}
\!+\! \left(b_{01}-b_{12}^db_{04},b_{02}, b_{03}\right) \cr
&+ \left( b_{12}-\alpha b_{13},b_{11}-\beta b_{13},
 b_{11}-b_{14},
 b_{2i}-\alpha
\right). \cr
}
$$
Thus we have proved:

\thm
\label{\thmthree}
For $n \ge 3$,
$$
\eqalignno{
\Ass &\left( {R \over K(n,d)} \right) \subseteq
\{Q_{i\Lambda}, Q_j, Q_{k\Lambda'}|
i = 4, 6, 10,11,12,13,14; j = 2,3,5; k = 1, 7, 9\} \cr
&\cup
\{Q_{8\Lambda\alpha},Q_{15\Lambda'\alpha'},Q_{16\Lambda'\alpha}\}
\cup
\lbrace Q_{jt}, Q_{19t\alpha\beta} | j = 17, 18;t = 2 \ldots, n \rbrace
\cr
&\cup \Ass \left({R \over K(n-1,d^2) + \C_1 + (b_{01}, b_{02}, b_{03}, b_{04})} \right), \cr
}
$$
where $\alpha$ and $\beta$ vary over all the $d$th roots of unity in $K$,
$\alpha'$ varies over all the $d^2$ roots of unity in $K$
which are not $d$th roots of unity,
$\Lambda$ varies over all the subsets of $\{1,2,3,4\}$,
and $\Lambda'$ varies over all the non-empty subsets of $\{1,2,3,4\}$.
\qed
\endb

It remains to find the associated primes of $K(n,d)$ when $n = 2$.
By Proposition~\propIrr,
it remains to find the associated primes of $K(2,d) : b_{04}^d c_{12}$,
which was computed on page~\qKEFcolon:
$$
\eqalignno{
&=
 b_{12}^{d^2} \left(b_{04}^d c_{11}-b_{01}^d c_{12},
b_{01}^d (c_{13}-c_{12}),
(b_{12}-b_{n-1,i})b_{03},
(b_{12}^d -b_{1i}^d )b_{04}, c_{13}-c_{12}\right) \cr
&\hskip2em+ (c_{11}-b_{12}^{d^2}c_{12},
c_{14}-c_{11}, b_{12}-b_{13},
b_{02}-b_{12}b_{03},
c_{12}b_{12}^{d^2}(b_{12}-b_{11})) \cr
&\hskip2em +(c_{12}b_{12}^{d^2}(b_{12}-b_{14}),
b_{01}-b_{12}^d b_{04},
c_{12}b_{11}-c_{13}b_{14},
b_{11}b_{13}^{d^2}-b_{14}b_{12}^{d^2}) \cr
&=
b_{12}^{d^2}\left((b_{12}-b_{11})c_{12}, b_{11}-b_{14},
(b_{12}-b_{1i})b_{03},
(b_{12}^d -b_{1i}^d )b_{04}, c_{13}-c_{12}\right) \cr
&\hskip2em
+ (c_{11}-b_{12}^{d^2}c_{12},
c_{14}-c_{11}, b_{12}-b_{13}, b_{02}-b_{12}b_{03},
b_{01}-b_{12}^d b_{04},
c_{12}b_{11}-c_{13}b_{14}). \cr
}
$$
By Fact~\factassocses,
any associated prime of $K(2,d) : b_{04}^d c_{12}$ is associated either to
$(K(2,d) : b_{04}^d c_{12}) + (b_{12}^{d^2})$ or to $K(2,d) : b_{04}^d c_{12} b_{12}^{d^2}$.
The former ideal is
$$
\eqalignno{
&=(b_{12}^{d^2},c_{11},
c_{14}, b_{12}-b_{13}, b_{02}-b_{12}b_{03},
b_{01}-b_{12}^d b_{04},
c_{12}b_{11}-c_{13}b_{14}), \cr
}
$$
which is a primary ideal,
so this contributes only the following prime ideal
to a list of possibly associated primes of $K(2,d)$:
$$
\eqalignno{
Q_5&=(c_{11},c_{14},b_{12},
b_{13}, b_{01},b_{02}, c_{12}b_{11}-c_{13}b_{14}). \cr
}
$$
Next we analyze $K(2,d) : b_{04}^d c_{12} b_{12}^{d^2}$:
$$
\eqalignno{
K(2,d) &: b_{04}^d c_{12} b_{12}^{d^2}=
\left((b_{12}-b_{11})c_{12}, b_{11}-b_{14},(b_{12}-b_{1i})b_{03},
(b_{12}^d -b_{1i}^d )b_{04}, c_{13}-c_{12}\right) \cr
&\hskip2em+ (c_{11}-b_{12}^{d^2}c_{12},
c_{14}-c_{11}, b_{12}-b_{13}, b_{02}-b_{12}b_{03},
b_{01}-b_{12}^d b_{04},
c_{12}b_{11}-c_{13}b_{14}) \cr
&=
\left((b_{12}-b_{11})c_{12}, b_{11}-b_{14}(b_{12}-b_{1i})b_{03},
(b_{12}^d -b_{1i}^d )b_{04}, c_{13}-c_{12}\right)\cr
&\hskip2em+ (c_{11}-b_{12}^{d^2}c_{12},
c_{14}-c_{11}, b_{12}-b_{13}, b_{02}-b_{12}b_{03},b_{01}-b_{12}^d b_{04}). \cr
}
$$
This decomposes further as
(first coloning and adding $c_{12}$):
$$
\eqalignno{
&=
\Bigl(b_{12}-b_{1i}, c_{13}-c_{12},
c_{11}-b_{12}^{d^2}c_{12},
c_{14}-c_{11}, b_{02}-b_{12}b_{03},
b_{01}-b_{12}^d b_{04} \Bigr) \cr
&\hskip1em\cap \Bigl(\C_1 + (b_{11}-b_{14},
b_{12}-b_{13}, b_{02}-b_{12}b_{03},b_{01}-b_{12}^d b_{04},
(b_{12}-b_{11})b_{03},
(b_{12}^d -b_{11}^d )b_{04})\Bigr) \cr
&=
\Bigl(b_{12}-b_{1i}, c_{13}-c_{12},
c_{11}-b_{12}^{d^2}c_{12},
c_{14}-c_{11}, b_{02}-b_{12}b_{03},
b_{01}-b_{12}^d b_{04} \Bigr) \cr
&\hskip1em\cap \Bigl(\C_1 + (b_{11}-b_{1i},
b_{02}-b_{12}b_{03},
b_{01}-b_{12}^d b_{04})\Bigr) \cr
&\hskip1em\cap \Bigl(\C_1 + (b_{11}-b_{14},
b_{12}-b_{13}, b_{02},b_{03},(b_{01}-b_{12}^d b_{04},
(b_{12}^d -b_{11}^d )b_{04})\Bigr) \cr
&=
\Bigl(b_{12}-b_{1i}, c_{13}-c_{12},
c_{11}-b_{12}^{d^2}c_{12},
c_{14}-c_{11}, b_{02}-b_{12}b_{03},
b_{01}-b_{12}^d b_{04} \Bigr) \cr
&\hskip1em\cap \Bigl(\C_1 + (b_{11}-b_{1i},
b_{02}-b_{12}b_{03},
b_{01}-b_{12}^d b_{04})\Bigr) \cr
&\hskip1em\cap \Bigl(\C_1 + (b_{11}-b_{14},
b_{12}-b_{13}, b_{02},b_{03},b_{01}-b_{12}^d b_{04},
b_{12}^d -b_{11}^d)\Bigr) \cr
&\hskip1em\cap \Bigl(\C_1 + (b_{11}-b_{14},
b_{12}-b_{13}, b_{01},b_{02},b_{03},b_{04})\Bigr). \cr
}
$$
From this decomposition it is easy to read off the associated primes,
which are possibly associated to $K(2,d)$:
$$
\eqalignno{
Q_{17,2} &= \Bigl(b_{12}-b_{1i}, c_{13}-c_{12},
c_{11}-b_{12}^{d^2}c_{12},
c_{14}-c_{11}, b_{02}-b_{12}b_{03},
b_{01}-b_{12}^{d}b_{04} \Bigr) \cr
Q_{18,2} &= \C_1 + (b_{11}-b_{1i},
b_{02}-b_{12}b_{03},
b_{01}-b_{12}^d b_{04}), \cr
Q_{19,2,1,\alpha} &= \C_1 + (b_{11}-b_{14},
b_{12}-b_{13}, b_{02},b_{03},b_{01}-b_{12}^d b_{04},
b_{12}-\alpha b_{11}), \alpha^d  = 1, \cr
Q_{20} &= \C_1 + (b_{11}-b_{14},
b_{12}-b_{13}, b_{01},b_{02},b_{03},b_{04}. \cr
}
$$

This completes finding an upper bound on the set of associated primes of $K(n,d)$,
for all $n \ge 2$:

\thm
\label{\thmtwo}
For $n = 2$,
$$
\Ass \left( {R \over K(n,d)} \right) \subseteq
\{Q_{1\Lambda'}, Q_j, Q_{i,2}, Q_{19,2,1,\alpha}, Q_{20}|
i = 17, 18; j = 2,3,5\},
$$
where $\alpha$ varies all the $d$ roots of unity in $K$,
and $\Lambda'$ varies over all non-empty subsets of $\{1,2,3,4\}$.
\qed
\endb

\vskip 4ex

\bigskip
\leftline{\bf References}
\bigskip

\bgroup
\font\eightrm=cmr8 \def\rm{\fam0\eightrm}
\font\eightit=cmti8 \def\it{\fam\itfam\eightit}
\font\eightbf=cmbx8 \def\bf{\fam\bffam\eightbf}
\font\eighttt=cmtt8 \def\tt{\fam\ttfam\eighttt}
\rm
\baselineskip=9.9pt
\parindent=3.6em

\item{[BS]}
D.\ Bayer and M.\ Stillman,
On the complexity of computing syzygies,
{\it J.\ Symbolic Comput.}, {\bf 6} (1988), 135-147.

\item{[D]}
M.\ Demazure,
Le th\'eor\`eme de complexit\'e de Mayr et Meyer,
{\it G\'eom\'etrie alg\'ebrique et applications, I
(La R\'abida, 1984)}, 35-58,
Travaux en Cours, 22, {\it Hermann, Paris}, 1987.

\item{[GS]}
D.\ Grayson and M.\ Stillman,
Macaulay2. 1996.
A system for computation in algebraic geometry and commutative algebra,
available via anonymous {\tt ftp} from {\tt math.uiuc.edu}.

\item{[GPS]}
G.-M.\ Greuel, G.\ Pfister and H.\ Sch\"onemann,
Singular. 1995.
A system for computation in algebraic geometry and singularity theory.
Available via anonymous {\tt ftp} from {\tt helios.mathematik.uni-kl.de}.

\item{[H]}
G.\ Herrmann,
Die Frage der endlich vielen Schritte in der Theorie der Polynomideale,
{\it Math.\ Ann.}, {\bf 95} (1926), 736-788.

\item{[K]}
J.\ Koh,
Ideals generated by quadrics exhibiting double exponential degrees,
{\it J.\ Algebra}, {\bf 200} (1998), 225-245.

\item{[MM]}
E.\ Mayr and A.\ Meyer,
The complexity of the word problems for commutative semigroups
and polynomial ideals,
{\it Adv.\ Math.}, {\bf 46} (1982), 305-329.

\item{[S1]}
I.\ Swanson,
The first Mayr-Meyer ideal,
preprint, 2001.

\item{[S2]}
I.\ Swanson,
The minimal components of the Mayr-Meyer ideals,
preprint, 2002.

\item{[S3]}
I.\ Swanson,
On the embedded primes of the Mayr-Meyer ideals,
preprint, 2002.

\egroup

\vskip 4ex
\noindent
{\sl
New Mexico State University - Department of Mathematical Sciences,
Las Cruces, New Mexico 88003-8001, USA.
E-mail: {\tt iswanson@nmsu.edu}.
}

\end